\documentclass{amsart}

\usepackage{amsmath,amssymb,amsthm,amscd}

\usepackage{nicefrac}

\usepackage{pinlabel}

\theoremstyle{definition}

\newtheorem*{rem}{Remark}

\newtheorem*{ack}{Acknowledgements}

\def\co{\colon\thinspace}

\newcommand{\rme}{\mathrm{e}}

\newcommand{\N}{\mathbb{N}}

\newcommand{\R}{\mathbb{R}}

\begin{document}

\author{Hansj\"org Geiges}
\address{Mathematisches Institut, Universit\"at zu K\"oln,
Weyertal 86--90, 50931 K\"oln, Germany}
\email{geiges@math.uni-koeln.de}

\title{Isotopies vis-\`a-vis level-preserving embeddings}

\date{}

\begin{abstract}
Various standard texts on differential topology maintain
that the level-preserving map defined by the
track of an isotopy of embeddings is itself an embedding.
This note describes a simple counterexample to this assertion.
\end{abstract}

\subjclass[2010]{57R40, 57R52}

\maketitle


\section{Introduction}
Let $Q,M$ be differential manifolds. An \emph{embedding}
of $Q$ in $M$ is a smooth map $\varphi\co Q\rightarrow M$ such that
$\varphi(Q)\subset M$ is a submanifold and $\varphi$ a diffeomorphism
onto its image. Equivalently, an embedding can be characterised
as an immersion (i.e.\ a smooth map of rank equal to
$\dim Q$) that maps $Q$ homeomorphically onto its
image~\cite[Theorem~1.3.1]{hirs76}.

An \emph{isotopy} of embeddings is a smooth map
\[ F\co Q\times [0,1]\longrightarrow M\]
with the property that for each $t\in[0,1]$ the map
\[\begin{array}{ccc}
Q & \longrightarrow & M\\
x & \longmapsto     & F(x,t)
\end{array}\]
is an embedding. The \emph{track} of this isotopy $F$
is the level-preserving map
\[\begin{array}{rccc}
G\co & Q\times [0,1] & \longrightarrow & M\times [0,1]\\
     & (x,t)         & \longmapsto     & \bigl(F(x,t),t\bigr).
\end{array}\]

It is not difficult to see that any level-preserving embedding
\[ G\co Q\times [0,1]\longrightarrow M\times[0,1],\]
that is, any embedding
$G$ satisfying $G\bigl(Q\times\{t\}\bigr)\subset M\times\{t\}$,
is the track of an
isotopy~\cite[Lemma II.4.2]{kosi93}. At least two of the
standard texts on differential topology maintain that, conversely,
the track of an isotopy is always an embedding~\cite[p.~178]{hirs76},
\cite[p.~34]{kosi93}.
This claim persists in the literature, at times presented as `easy to prove';
see \cite[p.~164]{shas11} or \cite[p.~207]{mukh15}, for instance.

The aim of this brief note is to exhibit a simple counterexample to
this assertion; notice that in any such counterexample the manifold
$Q$ needs to be non-compact. Thus, we are going to
show that \emph{the level-preserving map defined by the track
of an isotopy of embeddings is not, in general, an embedding.}

\begin{rem}
This observation is not new, see
\cite[Aufgabe 9.14]{brja73}, but apparently not nearly as well known
as it ought to be. 
\end{rem}

\section{The example}
In the example we take $Q=\R^+$, the positive real numbers, and
$M=\R^2$.
\subsection{A family of bump functions}
We first construct, in the standard way, a family of bump functions
$h_t\co \R\rightarrow [0,1]$, where $0<t\leq 1$,
with support in the interval $[\nicefrac{t}{2},
\nicefrac{3t}{2}]$.
Start with the smooth monotonically increasing function $f$ defined by
\[ f(s):=\begin{cases}
0           & \text{for $s\leq 0$},\\
\rme^{-1/s} & \text{for $s>0$}.
\end{cases} \]
Then the function $g_t\co\R\rightarrow [0,1]$ defined by
\[ g_t(s):=\frac{f(s)}{f(s)+f(\nicefrac{t}{4}-s)}\]
interpolates smoothly and monotonically between the value $0$ for $s\leq 0$
and the value $1$ for $s\geq t/4$. Finally, set
\[ h_t(s):= g_t\bigl(s-\nicefrac{t}{2}\bigr)\cdot
g_t(\nicefrac{3t}{2}-s\bigr).\]
This bump function is identically $0$
outside the interval $]\nicefrac{t}{2}, 
\nicefrac{3t}{2}[$, identically $1$ on the interval
$[\nicefrac{3t}{4},\nicefrac{5t}{4}]$, and it interpolates
monotonically in between. In particular, we have $h_t(t)=1$.

The key to constructing the desired isotopy is the
observation that the bump can be made to `disappear'
at $t=0$, provided we restrict the domain of definition
to the positive real numbers. In other words, if we take
$h_0\equiv 0$, then the function $(x,t)\mapsto h_t(x)$ will
be smooth on $\R^+\times[0,1]$.
\subsection{The isotopy}
\begin{figure}[h]
\labellist
\small\hair 2pt
\pinlabel $1$ [r] at 22 325
\pinlabel $t$ [b] at 181 -2
\pinlabel $3/2$ [b] at 468 -8
\pinlabel $2$ [b] at 612 -2
\endlabellist
\centering
\includegraphics[scale=0.4]{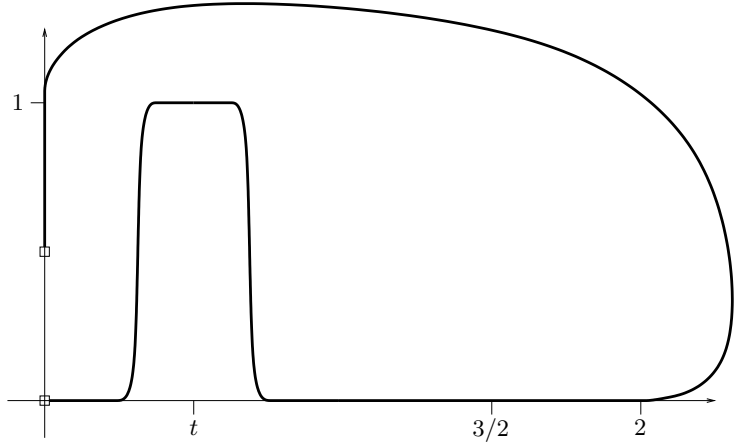}
  \caption{The image of the embedding $x\mapsto F(x,t)$, $0<t\leq 1$.}
  \label{figure:embedding}
\end{figure}

We now construct an isotopy $F\co\R^+\times [0,1]\rightarrow\R^2$
of embeddings $\R^+\rightarrow\R^2$. Set
\[ F(x,t):=\begin{cases}
\bigl(x,h_t(x)\bigr)       & \text{for $x\in]0,2]$ and $t\in]0,1]$},\\
(x,0)                      & \text{for $x\in]0,2]$ and $t=0$},\\
\text{independent of $t$}, & \\
\text{as shown in Figure~\ref{figure:embedding}} & \text{for $x\geq 2$}.
\end{cases}\]
Here it is understood that the choice for $x\geq 2$ is made in such a way
that $x\mapsto F(x,t)$ is an embedding.
The following three points, all of which are a simple
consequence of the definition, establish that $F$ is an isotopy.
\begin{itemize}
\item[(i)] The map $F$ is smooth at any point $(x,t)$
with $t>0$ or $x>\nicefrac{3}{2}$.
\item[(ii)] The map $F$ is also smooth at any given point $(x,0)$
with $0<x\leq\nicefrac{3}{2}$, since for any sequence $(x_{\nu},t_{\nu})$
converging to $(x,0)$ we have $F(x_{\nu},t_{\nu})=0$ for
$\nu$ sufficiently large.
\item[(iii)] The map $x\mapsto F(x,t)$ is an embedding
$\R^+\rightarrow\R^2$ for any $t\in[0,1]$.
\end{itemize}

However, the track $G\co (x,t)\mapsto \bigl(F(x,t),t\bigr)$ of this
isotopy is not an embedding, since it is not a homeomorphism
of $\R^+\times [0,1]$ onto the image of~$G$. To see this,
choose the unique $x_0>2$ such that $F(x_0,t)=(0,1)\in\R^2$
for any $t\in [0,1]$. Then $G(x_0,0)=(0,1,0)$. Now let $(x_{\nu})_{\nu\in\N}$
be a sequence of real numbers in $]0,1]$ converging to~$0$.
Since, by construction, we have $h_t(t)=1$ for $0<t\leq 1$, 
\[ G(x_{\nu},x_{\nu})=\bigl(F(x_{\nu},x_{\nu}),x_{\nu}\bigr)=
(x_{\nu},1,x_{\nu}),\]
which converges to $(0,1,0)=G(x_0,0)$. But the
sequence $(x_{\nu},x_{\nu})$ in $\R^+\times [0,1]$ does not
converge to $(x_0,0)$.

\begin{rem}
(1) The crucial feature of this example is that the embeddings
of $\R^+$ in $\R^2$ given by $x\mapsto F(x,t)$ are not proper
for $x$ near~$0$. It is this improperness
which allows one to make a bump `disappear' smoothly in finite time.
The properness or not for $x$ near $\infty$, on the other hand, is irrelevant.

(2) Another way to describe the essential characteristic of this
example is to observe that the map
\[ t\longmapsto\bigl\{x\longmapsto F(x,t)\bigr\}\]
from the interval $[0,1]$ into the space of embeddings $\R^+\rightarrow\R^2$
is not continuous at $t=0$ when the space of smooth
maps $\R^+\rightarrow\R^2$ is equipped with the \emph{strong} topology
in the sense of~\cite[Chapter~2.1]{hirs76}.
\end{rem}

The recent text by Wall~\cite{wall16},
largely based on notes from the 1960s, avoids the pitfall
described in this note
by working directly with a notion of \emph{diffeotopy
of $Q$ in~$M$}, defined as a level-preserving embedding
$Q\times [0,1]\rightarrow M\times[0,1]$. Remark~(2) serves to
indicate that for non-compact manifolds
this is in some sense a more appropriate definition.

\begin{ack}
I thank Peter Albers, Fran\c{c}ois Laudenbach, Manfred Lehn and
Thomas Rot for useful conversations and correspondence. This note
was prompted by Tho\-mas Rot pointing out a related phenomenon to me that
he observed during his ongoing work. A search of the
literature then revealed that the erroneous assertion about tracks
of isotopies continues to appear in text-books. So the
example in this note, which I have been discussing in various
courses on differential topology over the past 20 years,
is for the record.
\end{ack}

\end{document}